\documentclass[reqno]{amsart}
\usepackage{amsmath}
\usepackage{amsthm}
\usepackage{amsfonts}
\usepackage{amssymb}
\usepackage{latexsym}
\usepackage{amscd}
\usepackage{stmaryrd}

\theoremstyle{plain}
\newtheorem{thm}{Theorem}[section]
\newtheorem{cor}[thm]{Corollary}
\newtheorem{lem}[thm]{Lemma}
\newtheorem{prop}[thm]{Proposition}
\newtheorem{rem}[thm]{Remark}

\numberwithin{equation}{section}

\newfont{\scyr}{wncyr10 scaled 550}
\def\shuffle{\,\mbox{\bf \scyr X}\,}
\def\tshuffle{\overset{t}{\shuffle}}
\def\tast{\overset{t}{\ast}}

\def\proof{\noindent {\bf Proof.\;}}

\def\wt{\operatorname{wt}}
\def\dep{\operatorname{dep}}
\def\height{\operatorname{ht}}
\def\Li{\operatorname{Li}}

\allowdisplaybreaks

\begin{document}

\title{Some relations of interpolated multiple zeta values}

\date{\today\thanks{The first author is supported by the National Natural Science Foundation of
China (Grant No. 11471245) and the Natural Science Foundation of Shanghai (grant no. 14ZR1443500).
The authors thank the anonymous referee for his/her helpful comments, which improve the paper greatly. }}

\author{Zhonghua Li \quad and \quad Chen Qin}

\address{School of Mathematical Sciences, Tongji University, No. 1239 Siping Road,
Shanghai 200092, China}

\email{zhonghua\_li@tongji.edu.cn}

\address{School of Mathematical Sciences, Tongji University, No. 1239 Siping Road,
Shanghai 200092, China}

\email{2014chen\_qin@tongji.edu.cn}



\begin{abstract}
In this paper, the extended double shuffle relations for interpolated multiple zeta values are established. As an application, Hoffman's relations for interpolated multiple zeta values are proved. Furthermore, a generating function for sums of interpolated multiple zeta values of fixed weight, depth and height is represented by hypergeometric functions, and we discuss some special cases.
\end{abstract}

\maketitle

%%---------------------------------------------------------------------------
%%------------------------Content-------------------------------------------
%%----------------------------------------------------------------------------

\section{Introduction}\label{Sec:Intro}

Let $\mathbf{k}=(k_1,k_2,\ldots,k_n)$ be a sequence of positive integers with $n\geqslant 1$, we define its weight, depth and height respectively by
$$\wt(\mathbf{k})=k_1+k_2+\cdots+k_n,\quad \dep(\mathbf{k})=n,\quad \height(\mathbf{k})=\sharp\{i\mid 1\leqslant i\leqslant n,k_i\geqslant 2\}.$$
If $k_1\geqslant 2$, then $\mathbf{k}$ is said to be admissible. For such an admissible index $\mathbf{k}$, the multiple zeta value (MZV for short) and the multiple zeta-star value (MZSV for short) indexed by $\mathbf{k}$ are defined as
\begin{align*}
\zeta(\mathbf{k})=\zeta(k_1,k_2,\ldots,k_n)=\sum\limits_{m_1>m_2>\cdots>m_n>0}\frac{1}{m_1^{k_1}m_2^{k_2}\cdots m_n^{k_n}}
\end{align*}
and
\begin{align*}
\zeta^{\star}(\mathbf{k})=\zeta^{\star}(k_1,k_2,\ldots,k_n)=\sum\limits_{m_1\geqslant m_2\geqslant \cdots\geqslant m_n\geqslant 1}\frac{1}{m_1^{k_1}m_2^{k_2}\cdots m_n^{k_n}},
\end{align*}
respectively. Let $t$ be a variable. In \cite{Yamamoto}, S. Yamamoto introduced the following definition
\begin{align}
\zeta^t(\mathbf{k})=\zeta^t(k_1,k_2,\ldots,k_n)=\sum\limits_{\mathbf{p}} t^{n-\dep(\mathbf{p})}\zeta(\mathbf{p})(\in\mathbb{R}[t]),
\label{Eq:tMZV}
\end{align}
where $\mathbf{p}$ runs over all sequences of the form
$$\mathbf{p}=(k_1\Box k_2\Box\cdots\Box k_n)$$
in which each $\Box$ is filled by the comma, or the plus $+$. Since
$$\zeta^0(\mathbf{k})=\zeta(\mathbf{k}),\quad \zeta^1(\mathbf{k})=\zeta^{\star}(\mathbf{k}),$$
we can regard $\zeta^t(\mathbf{k})$ as an interpolation polynomial of MZVs and MZSVs. We call the polynomials defined by \eqref{Eq:tMZV} interpolated multiple zeta values ($t$-MZVs for short).

There are several relations of $t$-MZVs which have been found. For example, for any integers $k>n\geqslant 1$, the sum formula
\begin{align}
\sum\limits_{\wt(\mathbf{k})=k,\dep(\mathbf{k})=n\atop \mathbf{k}:\text{admissible}}\zeta^t(\mathbf{k})=\left(\sum\limits_{i=0}^{n-1}\binom{k-1}{i}t^i(1-t)^{n-1-i}\right)\zeta(k)
\label{Eq:Sum-Formula}
\end{align}
was proved in \cite{Yamamoto}. Taking $t=0$, one obtains the sum formula for MZVs which was first proved in \cite{Granville}; and taking $t=1$, one gets sum formula for MZSVs which was first proved in \cite{Kombu}. In \cite{Yamamoto}, S. Yamamoto obtained the cyclic sum formula for $t$-MZVs:
\begin{align*}
&\sum\limits_{i=1}^n\sum\limits_{j=1}^{k_i-1}\zeta^t(k_i+1-j,k_{i+1},\ldots,k_n,k_1,\ldots,k_{i-1},j)\\
=&(1-t)\sum\limits_{i=1}^n\zeta^t(k_i+1,k_{i+1},\ldots,k_n,k_1,\ldots,k_{i-1})+t^nk\zeta(k+1),
\end{align*}
where $k_1,\ldots,k_n$ are positive integers with at least one greater than $1$, and $k=k_1+\cdots+k_n$. The cyclic sum formulas for MZVs and MZSVs were first proved separately in \cite{Hoffma-Ohno} and \cite{Ohno-Wakabayashi}. In \cite{Tanaka-Wakabayashi}, T. Tanaka and N. Wakabayashi established Kawashima's relations for $t$-MZVs and as an application, they gave a new proof of the cyclic sum formula. Kawashima's relations for MZVs were first appeared in \cite{Kawashima}, and the equivalent formulas for MZSVs were given in \cite{Tanaka}.

We give some other relations among $t$-MZVs in this paper. In Section \ref{Sec:EDS}, we discuss the extended double shuffle relations. The extended double shuffle relations for MZVs were developed  in \cite{Ihara-Kaneko-Zagier,Racinet}, and for MZSVs were introduced in \cite{Muneta}. By a similar way as in \cite{Muneta}, and with the helps of \cite{Tanaka-Wakabayashi,Yamamoto}, we get the extended double shuffle relations (Theorem \ref{Thm:EDS}) for $t$-MZVs. As an application of our extended double shuffle relations, we also get Hoffman's relations for $t$-MZVs.

Then we discuss the sums of $t$-MZVs of fixed weight, depth and height in Section \ref{Sec:Sum} by a similar method as in \cite{Ohno-Zagier,Aoki-Kombu-Ohno,Aoki-Ohno}. There are some related work  on such type of sums. In \cite{Ohno-Zagier}, Y. Ohno and D. Zagier showed that the sum of MZVs of fixed weight, depth and height can be expressed as a polynomial of Riemann zeta values with rational coefficients. In \cite{Aoki-Kombu-Ohno,Aoki-Ohno}, T. Aoki, Y. Kombu and Y. Ohno expressed a generating function of sums of MZSVs of fixed weight, depth and height by hypergeometric functions and obtained some relations of MZSVs.  We discussed sum of MZVs of fixed weight, depth and $i$-height in \cite{Li2008}, and obtained its $q$-analogue in \cite{Li2010}. And T. Aoki, Y. Ohno and N. Wakabayashi studied the sums of MZSVs of fixed weight, depth and $i$-height in \cite{Aoki-Ohno-Wakabayashi}.

Now for positive integers $k,n,s$ with $k\geqslant n+s$ and $n\geqslant s$, we define a sum
$$X_0(k,n,s)=\sum\limits_{\mathbf{k}\in I_0(k,n,s)}\zeta^t(\mathbf{k}),$$
where $I_0(k,n,s)$ is the set of all admissible indices with weight $k$, depth $n$ and height $s$. Then a generating function of these sums is defined as
$$\Phi_0(u,v,w)=\sum\limits_{k\geqslant n+s,n\geqslant s\geqslant 1}X_0(k,n,s)u^{k-n-s}v^{n-s}w^{2s-2},$$
where $u,v,w$ are variables. We express this generating function by hypergeometric functions (Theorem \ref{Thm:Sum}) and discuss some special cases.  In the case $uv=w^2$, we give a new proof of the sum formula \eqref{Eq:Sum-Formula}. In the case $v=0$, we find that the sums $X_0(k,n,n)$ can be expressed as  polynomials of Riemann zeta values with $\mathbb{Q}[t]$ coefficients. And in the case $w=0$, we give an expression for height one $t$-MZVs.

There are two appendices. In Appendix \ref{AppSec:Proof-Sum}, we give another proof of the sum formula \eqref{Eq:Sum-Formula} by the sum formula of MZVs. And in Appendix \ref{AppSec:Proof-Gauss}, we give a proof of \eqref{Eq:Gauss-Hyper-Duality}, which is an identity of Gauss hypergeometric functions.

\begin{rem}
After the first version of this paper was submitted to arXiv, N. Wakabayashi kindly informed us her work \cite{Wakabayashi}, in which she independently got the extended double shuffle and Hoffman's relations of $t$-MZVs.
\end{rem}

%%%%%%%%%%%%%%%%%%%%%%%%%%%%%%%%%%%%%%%%%%%%%%%%%%%%%%%%%%%%%%%%

\section{Extended double shuffle relations of $t$-MZVs} \label{Sec:EDS}

In this section, we discuss the extended double shuffle relations of $t$-MZVs.

\subsection{Algebraic setup}

As in \cite{Hoffman1997,Hoffma-Ohno,Ihara-Kaneko-Zagier,Muneta,Tanaka-Wakabayashi,Yamamoto}, we give the algebraic setup. Let $A=\{x,y\}$ be an alphabet with two noncommutative letters, and denote by $A^{\ast}$ the set of all words generated by $A$, which contains the empty word $1$. Let $\mathfrak{h}_t=\mathbb{Q}[t]\langle A\rangle$ be the noncommutative polynomial algebra over $\mathbb{Q}[t]$ generated by $A$, and define two subalgebras
$$\mathfrak{h}_t^1=\mathbb{Q}[t]+\mathfrak{h}_ty,\qquad \mathfrak{h}_t^{0}=\mathbb{Q}[t]+x\mathfrak{h}_ty.$$
For any positive integer $k$, set $z_k=x^{k-1}y\in A^{\ast}$. We define a $\mathbb{Q}[t]$-linear map $Z_t:\mathfrak{h}^{0}_t\rightarrow\mathbb{R}[t]$ by $Z_t(1)=1$ and
$$Z_t(z_{k_1}\cdots z_{k_n})=\zeta^t(k_1,\ldots,k_n),\quad (n,k_1,\ldots,k_n\in\mathbb{N},k_1\geqslant 2).$$

We denote $\mathfrak{h}_0$, $\mathfrak{h}_0^1$, $\mathfrak{h}_0^0$ and $Z_0$ simply by $\mathfrak{h}$, $\mathfrak{h}^1$, $\mathfrak{h}^0$ and $Z$, respectively. Hence it is obvious that
$$\mathfrak{h}_t=\mathfrak{h}[t],\quad \mathfrak{h}_t^1=\mathfrak{h}^1[t],\quad \mathfrak{h}_t^0=\mathfrak{h}^0[t].$$
Let $\sigma_t$ be the automorphism of the algebra $\mathfrak{h}_t$ determined by
$$\sigma_t(x)=x,\quad \sigma_t(y)=tx+y.$$
Note that $\sigma_t^{-1}=\sigma_{-t}$. We now define a $\mathbb{Q}[t]$-linear map $S_t:\mathfrak{h}_t\rightarrow \mathfrak{h}_t$ by $S_t(1)=1$ and
$$S_t(wa)=\sigma_t(w)a,\quad (w\in \mathfrak{h}_t,a\in A).$$
Note that the restriction of $S_t$ on $\mathfrak{h}_t^1$ is nothing but that defined in \cite{Tanaka-Wakabayashi,Yamamoto}. Then it is easy to prove that $S_t$ is invertible with the inverse $S_t^{-1}=S_{-t}$, and
$$S_t(\mathfrak{h}_t^1)=\mathfrak{h}_t^1,\qquad S_t(\mathfrak{h}_t^0)=\mathfrak{h}_t^0.$$
We extend the $\mathbb{Q}$-linear map $Z:\mathfrak{h}^0\rightarrow \mathbb{R}$ to a $\mathbb{Q}[t]$-linear map $Z:\mathfrak{h}^0_t=\mathfrak{h}^0[t]\rightarrow \mathbb{R}[t]$, then one can show that on $\mathfrak{h}_t^0$
$$Z_t=Z\circ S_t.$$

We define a new product $\tshuffle$ on the space $\mathfrak{h}_t$, which we call $t$-shuffle product. The definition here is similar to that one given in \cite{Muneta} for MZSVs, but with some modifications. The $t$-shuffle product $\tshuffle:\mathfrak{h}_t\times \mathfrak{h}_t\rightarrow\mathfrak{h}_t$ is $\mathbb{Q}[t]$-bilinear, and satisfies the rules
\begin{itemize}
  \item[(S1)] $1\tshuffle w=w\tshuffle 1=w$,
  \item[(S2)] $aw_1\tshuffle bw_2=a(w_1\tshuffle bw_2)+b(aw_1\tshuffle w_2)-\delta(w_1)\rho(a)bw_2-\delta(w_2)\rho(b)aw_1$,
\end{itemize}
where $w,w_1,w_2\in A^{\ast}$, $a,b\in A$, the map $\delta:A^{\ast}\rightarrow \{0,1\}$ is defined by
$$\delta(w)=\begin{cases}
1 & \text{if\;} w=1,\\
0 & \text{if\;} w\neq 1,
\end{cases}$$
and the map $\rho:A\rightarrow \mathfrak{h}_t$ is defined by
$$\rho(x)=0,\quad \rho(y)=tx.$$
Note that in the case $t=1$, the map $\rho$ is different from that of \cite{Muneta}, where the notation $\tau$ was used, and $\tau(x)=y$. By the definitions, it is easy to check that for any $a\in A$, we have
\begin{align}
\sigma_t(a)-\sigma_t(\rho(a))=a.
\label{Eq:sigma-rho}
\end{align}

When $t=0$, we get the usual shuffle product $\shuffle=\overset{0}{\shuffle}$ defined on the space $\mathfrak{h}$, which is commutative and associative. We can regard $\shuffle$ as a $\mathbb{Q}[t]$-bilinear product on the space $\mathfrak{h}_t$. Then we have

\begin{prop}\label{Prop:t-shuffle}
For any $w_1,w_2\in\mathfrak{h}_t$, we have
$$S_t(w_1\tshuffle w_2)=S_t(w_1)\shuffle S_t(w_2).$$
\end{prop}

\proof We may assume that $w_1,w_2\in A^{\ast}$. If $w_1=1$ or $w_2=1$, it is easy to get the result. Now assume that $w_1=u_1b_1$, $w_2=u_2b_2$ with $b_1,b_2\in A$ and $u_1,u_2\in A^{\ast}$. We use induction on $|u_1|+|u_2|$. Here for a word $w\in A^{\ast}$, we denote by $|w|$ the number of letters contained in $w$.

If $u_1=u_2=1$, the result follows from
\begin{align*}
S_t(x\tshuffle x)=&S_t(2x^2)=2x^2=x\shuffle x,\\
S_t(x\tshuffle y)=&S_t(y\tshuffle x)=S_t(xy+yx-tx^2)=xy+(tx+y)x-tx^2\\
=&xy+yx=x\shuffle y=y\shuffle x,\\
S_t(y\tshuffle y)=&S_t(2y^2-2txy)=2(tx+y)y-2txy=2y^2=y\shuffle y.
\end{align*}
If $u_1=1$ and $u_2=aw$ with $a\in A$ and $w\in A^{\ast}$, we have
\begin{align*}
&S_t(b_1\tshuffle awb_2)=S_t\left(b_1awb_2+a(b_1\tshuffle wb_2)-\rho(b_1)awb_2\right)\\
=&\sigma_t(b_1)S_t(awb_2)+\sigma_t(a)S_t(b_1\tshuffle wb_2)-\sigma_t\left(\rho(b_1)\right)S_t(awb_2).
\end{align*}
Using the induction hypothesis and \eqref{Eq:sigma-rho}, we get
$$S_t(b_1\tshuffle awb_2)=b_1S_t(awb_2)+\sigma_t(a)\left(b_1\shuffle S_t(wb_2)\right).$$
While we have
$$b_1\shuffle S_t(awb_2)=b_1\shuffle \sigma_t(a)S_t(wb_2)=b_1\sigma_t(a)S_t(wb_2)+\sigma_t(a)(b_1\shuffle S_t(wb_2)).$$
Hence we get the result in this case. Similarly, one can prove that the result holds for the case $|u_1|>0$ and $u_2=1$. Finally, if $u_1=a_1v_1$ and $u_2=a_2v_2$ with $a_1,a_2\in A$ and $v_1,v_2\in A^{\ast}$, we have
\begin{align*}
&S_t(a_1v_1b_1\tshuffle a_2v_2b_2)=S_t\left(a_1(v_1b_1\tshuffle a_2v_2b_2)+a_2(a_1v_1b_1\tshuffle v_2b_2)\right)\\
=&\sigma_t(a_1)S_t(v_1b_1\tshuffle a_2v_2b_2)+\sigma_t(a_2)S_t(a_1v_1b_1\tshuffle v_2b_2).
\end{align*}
By the induction hypothesis, we get
\begin{align*}
&S_t(a_1v_1b_1\tshuffle a_2v_2b_2)\\
=&\sigma_t(a_1)(S_t(v_1b_1)\shuffle S_t(a_2v_2b_2))+\sigma_t(a_2)(S_t(a_1v_1b_1)\shuffle S_t(v_2b_2))\\
=&S_t(a_1v_1b_1)\shuffle S_t(a_2v_2b_2),
\end{align*}
which finishes the proof.
\qed

Since the shuffle product $\shuffle$ is commutative and associative, and the map $Z:(\mathfrak{h}^0,\shuffle)\rightarrow \mathbb{R}$ is an algebra homomorphism, one can get the following theorem from Proposition \ref{Prop:t-shuffle} without difficulty.

\begin{thm}\label{Thm:t-shuffle-Alg}
The $t$-shuffle product $\tshuffle$ is commutative and associative, $\mathfrak{h}_t$ is a commutative and associative algebra under this product, $\mathfrak{h}_t^1$ and $\mathfrak{h}^0_t$ are still subalgebras. Moreover, the $\mathbb{Q}[t]$-linear map $Z_t:(\mathfrak{h}_t^0,\tshuffle)\rightarrow\mathbb{R}[t]$ is an algebra homomorphism.
\end{thm}

As an application of above discussions, we derive Euler's decomposition formula for $t$-MZVs.

\begin{prop}\label{Prop:Euler-Decom}
For any positive integers $k$ and $l$, we have
\begin{align*}
z_k\tshuffle z_l=&\sum\limits_{i=1}^{k}\binom{k+l-i-1}{l-1}z_{k+l-i}z_i\\
&\quad+\sum\limits_{i=1}^l\binom{k+l-i-1}{k-1}z_{k+l-i}z_i-\binom{k+l}{k}tz_{k+l}.
\end{align*}
Moreover, when $k,l\geqslant 2$, we have
\begin{align*}
\zeta(k)\zeta(l)=&\sum\limits_{i=1}^{k}\binom{k+l-i-1}{l-1}\zeta^t(k+l-i,i)\\
&\quad+\sum\limits_{i=1}^l\binom{k+l-i-1}{k-1}\zeta^t(k+l-i,i)-\binom{k+l}{k}t\zeta(k+l).
\end{align*}
\end{prop}

\proof The usual Euler's decomposition formula is
$$z_k\shuffle z_l=\sum\limits_{i=1}^{k}\binom{k+l-i-1}{l-1}z_{k+l-i}z_i+\sum\limits_{i=1}^l\binom{k+l-i-1}{k-1}z_{k+l-i}z_i.$$
Hence we have
\begin{align*}
&z_k\tshuffle z_l=S_{-t}(S_t(z_k)\shuffle S_t(z_l))=S_{-t}(z_k\shuffle z_l)\\
=&\sum\limits_{i=1}^{k}\binom{k+l-i-1}{l-1}S_{-t}(z_{k+l-i}z_i)+\sum\limits_{i=1}^l\binom{k+l-i-1}{k-1}S_{-t}(z_{k+l-i}z_i)\\
=&\sum\limits_{i=1}^{k}\binom{k+l-i-1}{l-1}z_{k+l-i}z_i+\sum\limits_{i=1}^l\binom{k+l-i-1}{k-1}z_{k+l-i}z_i\\
&-t\left(\sum\limits_{i=1}^{k}\binom{k+l-i-1}{l-1}+\sum\limits_{i=1}^l\binom{k+l-i-1}{k-1}\right)z_{k+l},
\end{align*}
from which one easily gets the result.
\qed

Now we recall the definition and properties of $t$-harmonic shuffle product on $\mathfrak{h}_t^1$ from \cite{Tanaka-Wakabayashi,Yamamoto}.
The $\mathbb{Q}[t]$-bilinear product $\tast:\mathfrak{h}_t^1\times \mathfrak{h}_t^1\rightarrow\mathfrak{h}_t^1$ is defined by the rules
\begin{itemize}
  \item[(H1)] $1\tast w=w\tast 1=w$,
  \item[(H2)] $z_kw_1\tast z_lw_2=z_k(w_1\tast z_lw_2)+z_l(z_kw_1\tast w_2)+(1-2t)z_{k+l}(w_1\tast w_2)+(1-\delta(w_1)\delta(w_2))(t^2-t)x^{k+l}(w_1\tast w_2)$,
\end{itemize}
where $w,w_1,w_2\in A^{\ast}\cap\mathfrak{h}_t^1$ and $k,l\in \mathbb{N}$. When $t=0$, we get the usual harmonic shuffle product $\ast=\overset{0}{\ast}$ defined on the space $\mathfrak{h}^1$. We can regard $\ast$ as a $\mathbb{Q}[t]$-bilinear product on the space $\mathfrak{h}_t^1$.

\begin{prop}[\cite{Yamamoto}]\label{Prop:t-ast}
For any $w_1,w_2\in\mathfrak{h}_t^1$, we have
$$S_t(w_1\tast w_2)=S_t(w_1)\ast S_t(w_2).$$
\end{prop}

Since the harmonic shuffle product $\ast$ is commutative and associative, and the map $Z:(\mathfrak{h}^0,\ast)\rightarrow \mathbb{R}$ is an algebra homomorphism, the following theorem follows from Proposition \ref{Prop:t-ast}.

\begin{thm}[\cite{Tanaka-Wakabayashi,Yamamoto}]\label{Thm:t-ast-Alg}
The $t$-harmonic shuffle product $\tast$ is commutative and associative, $\mathfrak{h}_t^1$ is a commutative and associative algebra under this product, and $\mathfrak{h}^0_t$ is still a subalgebra. Moreover, the $\mathbb{Q}[t]$-linear map $Z_t:(\mathfrak{h}_t^0,\tast)\rightarrow\mathbb{R}[t]$ is an algebra homomorphism.
\end{thm}

Finally, By Theorem \ref{Thm:t-shuffle-Alg} and Theorem \ref{Thm:t-ast-Alg}, we get the finite double shuffle relations for $t$-MZVs.

\begin{thm}\label{Thm:FDS}
For any $w_1,w_2\in\mathfrak{h}_t^0$, we have
$$Z_t(w_1\tshuffle w_2-w_1\tast w_2)=0.$$
\end{thm}

As an example, by Euler's decomposition formula (Proposition \ref{Prop:Euler-Decom}) and the formula
$$z_t\tast z_l=z_kz_l+z_lz_k+(1-2t)z_{k+l},$$
we get the finite double shuffle relation
\begin{align*}
&\left\{1+\left[\binom{k+l}{k}-2\right]t\right\}\zeta(k+l)\\
=&\sum\limits_{i=1}^{k-1}\binom{k+l-i-1}{l-1}\zeta^t(k+l-i,i)+\sum\limits_{i=1}^{l-1}\binom{k+l-i-1}{k-1}\zeta^t(k+l-i,i),
\end{align*}
which holds for any integers $k,l\geqslant 2$.

\subsection{Extended double shuffle relations}

Similarly as in \cite{Muneta}, we can get extended double shuffle relations of $t$-MZVs from that of MZVs. We recall the extended double shuffle relations of MZVs (see \cite{Ihara-Kaneko-Zagier} for example). As shuffle algebras, we have $\mathfrak{h}^1=\mathfrak{h}^0[y]$. Therefore there is a unique algebra homomorphism $Z^{\shuffle}:(\mathfrak{h}^1,\shuffle)\rightarrow\mathbb{R}[T]$, such that
$$Z^{\shuffle}|_{\mathfrak{h}^0}=Z,\quad Z^{\shuffle}(y)=T.$$
Similarly, there is a unique algebra homomorphism $Z^{\ast}:(\mathfrak{h}^1,\ast)\rightarrow\mathbb{R}[T]$, such that
$$Z^{\ast}|_{\mathfrak{h}^0}=Z,\quad Z^{\ast}(y)=T.$$
Then the extended double shuffle relations for MZVs claim that for any $w_1\in\mathfrak{h}^1$ and $w_0\in\mathfrak{h}^0$, it holds
$$Z^{\shuffle}(w_1\shuffle w_0-w_1\ast w_0)=0,\quad\text{and}\quad Z^{\ast}(w_1\shuffle w_0-w_1\ast w_0)=0.$$

Now we return to $t$-MZVs. We have the following result, which is the $t$-MZVs counterpart of \cite[Lemma 2.10]{Muneta}.

\begin{lem}\label{Lem:Stru-1-0}
As $t$-shuffle algebras, we have $\mathfrak{h}_t^1=\mathfrak{h}_t^0[y]$. Similarly, as $t$-harmonic shuffle algebras, we have $\mathfrak{h}_t^1=\mathfrak{h}_t^0[y]$.
\end{lem}

\proof We treat $t$-shuffle algebras. It is enough to show that for any integer $m\geqslant 0$ and any $w\in\mathfrak{h}_t^0\cap A^{\ast}$, the word $y^mw$ can be uniquely expressed as
$$y^mw=\sum\limits_{i=0}^mw_i\tshuffle y^{\tshuffle i}$$
with $w_0,\ldots,w_m\in\mathfrak{h}_t^0$.

Since $S_t(y^mw)\in\mathfrak{h}_t^1$, by the fact $(\mathfrak{h}^1,\shuffle)=(\mathfrak{h}^0,\shuffle)[y]$ (more precisely \cite[Corollary 5]{Ihara-Kaneko-Zagier}), we have
$$S_t(y^mw)=\sum\limits_{i=0}^mu_i\shuffle y^{\shuffle i},\quad (\exists u_i\in\mathfrak{h}_t^0).$$
By Proposition \ref{Prop:t-shuffle}, we get
$$y^mw=\sum\limits_{i=0}^mS_{-t}(u_i\shuffle y^{\shuffle i})=\sum\limits_{i=0}^mS_{-t}(u_i)\tshuffle y^{\tshuffle i}.$$
Since $S_{-t}(u_i)\in\mathfrak{h}_t^0$, we get the existence of the representation.

For the uniqueness of the representation, assume that
$$\sum\limits_{i=0}^mw_i\tshuffle y^{\tshuffle i}=\sum\limits_{j=0}^nv_j\tshuffle y^{\tshuffle j},\quad (w_i,v_j\in\mathfrak{h}_t^0).$$
Applying the map $S_t$, we get
$$\sum\limits_{i=0}^mS_t(w_i)\shuffle y^{\shuffle i}=\sum\limits_{j=0}^nS_t(v_j)\shuffle y^{\shuffle j}.$$
Since $S_t(w_i),S_t(v_j)\in\mathfrak{h}_t^0$, we get $m=n$ and $S_t(w_i)=S_t(v_i)$ for all $i$. Hence $w_i=v_i$ for all $i$.
\qed

By Lemma \ref{Lem:Stru-1-0}, there exist unique algebra homomorphisms $Z_t^{\shuffle}:(\mathfrak{h}_t^1,\tshuffle)\rightarrow\mathbb{R}[t,T]$ and $Z_t^{\ast}:(\mathfrak{h}_t^1,\tast)\rightarrow \mathbb{R}[t,T]$, such that
$$Z_t^{\shuffle}|_{\mathfrak{h}_t^0}=Z_t=Z_t^{\ast}|_{\mathfrak{h}_t^0},\qquad Z_t^{\shuffle}(y)=Z_t^{\ast}(y)=T.$$
Note that it is easy to check that
$$Z_t^{\shuffle}=Z^{\shuffle}\circ S_t,\quad Z_t^{\ast}=Z^{\ast}\circ S_t.$$
Finally we obtain the extended double shuffle relations for $t$-MZVs.

\begin{thm}\label{Thm:EDS}
For any $w_1\in\mathfrak{h}^1_t$ and $w_0\in\mathfrak{h}^0_t$, we have
$$Z^{\shuffle}_t(w_1\tshuffle w_0-w_1\tast w_0)=0,\quad\text{and}\quad Z^{\ast}_t(w_1\tshuffle w_0-w_1\tast w_0)=0.$$
\end{thm}

\proof We prove the first identity. The proof of the second identity is similar. We have
\begin{align*}
&Z^{\shuffle}_t(w_1\tshuffle w_0-w_1\tast w_0)=Z^{\shuffle}\left(S_t(w_1\tshuffle w_0-w_1\tast w_0)\right)\\
=&Z^{\shuffle}(S_t(w_1)\shuffle S_t(w_0)-S_t(w_1)\ast S_t(w_0))=0,
\end{align*}
which uses the facts $S_t(\mathfrak{h}_t^1)=\mathfrak{h}_t^1$ and $S_t(\mathfrak{h}_t^0)=\mathfrak{h}_t^0$.
\qed

\subsection{Hoffman's relations}

As an application, we take $w_1=y$ in Theorem \ref{Thm:EDS}, and get the following result, which can be regraded as Hoffman's relations for $t$-MZVs.

\begin{thm}\label{Thm:Hoffman}
For any positive integers $n,k_1,\ldots,k_n$ with $k_1\geqslant 2$, we have
\begin{align*}
&\sum\limits_{i=1}^n\left[1+(k_i+\delta_{ni}-2)t\right]\zeta^t(k_1,\ldots,k_{i-1},k_i+1,k_{i+1},\ldots,k_n)\\
=&\sum\limits_{i=1}^n\sum\limits_{j=2}^{k_i}\zeta^t(k_1,\ldots,k_{i-1},j,k_{i}+1-j,k_{i+1},\ldots,k_n)\\
&+(t-t^2)\sum\limits_{i=1}^{n-1}\zeta^t(k_1,\ldots,k_{i-1},k_i+k_{i+1}+1,k_{i+2},\ldots,k_n),
\end{align*}
where $\delta_{ij}$ is the Kronecker delta symbol.
\end{thm}

Setting $t=0$ and $t=1$ in Theorem \ref{Thm:Hoffman}, one obtains Hoffman's relations for MZVs which were first proved in \cite{Hoffman1992} and Hoffman's relations for MZSVs which were first proved in \cite{Muneta}, respectively. Theorem \ref{Thm:Hoffman} is immediately deduced from Theorem \ref{Thm:EDS} and the following lemma.

\begin{lem}\label{Lem:y-z1-zn}
For positive integers $n,k_1,\ldots,k_n$, we have
\begin{align}
y\tshuffle z_{k_1}\cdots z_{k_n}=&\sum\limits_{i=0}^nz_{k_1}\cdots z_{k_i}z_1z_{k_{i+1}}\cdots z_{k_n}\nonumber\\
&+\sum\limits_{i=1}^n\sum\limits_{j=2}^{k_i}z_{k_1}\cdots z_{k_{i-1}}z_jz_{k_i+1-j}z_{k_{i+1}}\cdots z_{k_n}\nonumber\\
&-\sum\limits_{i=1}^{n}(k_i+\delta_{ni})tz_{k_1}\cdots z_{k_{i-1}}z_{k_i+1}z_{k_{i+1}}\cdots z_{k_n},
\label{Eq:y-z1-zn-shuffle}
\end{align}
and
\begin{align}
y\tast z_{k_1}\cdots z_{k_n}=&\sum\limits_{i=0}^nz_{k_1}\cdots z_{k_i}z_1z_{k_{i+1}}\cdots z_{k_n}\nonumber\\
&+(1-2t)\sum\limits_{i=1}^nz_{k_1}\cdots z_{k_{i-1}}z_{k_i+1}z_{k_{i+1}}\cdots z_{k_n}\nonumber\\
&+(t^2-t)\sum\limits_{i=1}^{n-1}z_{k_1}\cdots z_{k_{i-1}}z_{k_i+k_{i+1}+1}z_{k_{i+2}}\cdots z_{k_n}.
\label{Eq:y-z1-zn-ast}
\end{align}
\end{lem}

To prove \eqref{Eq:y-z1-zn-shuffle}, we need the following lemma.

\begin{lem}\label{Lem:y-tshuffle}
For any positive integer $k$ and any $w\in A^{\ast}$, we have
$$y\tshuffle z_kw=\begin{cases}
\sum\limits_{i=1}^kz_iz_{k+1-i}-(k+1)tz_{k+1}+z_kz_1 & \text{if\;} w=1,\\
\sum\limits_{i=1}^kz_iz_{k+1-i}w-ktz_{k+1}w+z_k(y\tshuffle w) & \text{if\;} w\neq 1.
\end{cases}$$
\end{lem}

\proof We prove the lemma by induction on $k$. For $k=1$, it is direct computations. Now assume that $k>1$, we have
$$y\tshuffle z_kw=y\tshuffle xz_{k-1}w=yz_kw+x(y\tshuffle z_{k-1}w)-tx^2z_{k-1}w.$$
Using induction hypothesis, we get
\begin{align*}
y\tshuffle z_k=&z_1z_k+\sum\limits_{i=1}^{k-1}xz_iz_{k-i}-ktxz_{k}+xz_{k-1}z_1-tz_{k+1}\\
=&\sum\limits_{i=1}^kz_iz_{k+1-i}-(k+1)tz_{k+1}+z_kz_1.
\end{align*}
Similarly, using induction hypothesis, one can prove the result for $w\neq 1$.
\qed

Now we give a proof of Lemma \ref{Lem:y-z1-zn}.

\noindent {\bf Proof of Lemma \ref{Lem:y-z1-zn}.} To prove \eqref{Eq:y-z1-zn-shuffle}, we use induction on $n$. If $n=1$, we get \eqref{Eq:y-z1-zn-shuffle} from the formula $y\tshuffle z_k$ given in Lemma \ref{Lem:y-tshuffle}. Now Assume that $n>1$. By Lemma \ref{Lem:y-tshuffle}, we have
\begin{align*}
y\tshuffle z_{k_1}\cdots z_{k_n}=&\sum\limits_{j=1}^{k_1}z_jz_{k_1+1-j}z_{k_2}\cdots z_{k_n}-k_1tz_{k_1+1}z_{k_2}\cdots z_{k_n}\\
&\quad+z_{k_1}(y\tshuffle z_{k_2}\cdots z_{k_n}).
\end{align*}
Then it is easy to get \eqref{Eq:y-z1-zn-shuffle} by induction hypothesis.

Similarly, we prove \eqref{Eq:y-z1-zn-ast} by induction on $n$. If $n=1$, we get \eqref{Eq:y-z1-zn-ast} by the definition of $t$-harmonic shuffle product. If $n>1$, we have
\begin{align*}
y\tast z_{k_1}\cdots z_{k_n}=&z_1z_{k_1}\cdots z_{k_n}+z_{k_1}(y\tast z_{k_2}\cdots z_{k_n})+(1-2t)z_{k_1+1}z_{k_2}\cdots z_{k_n}\\
&\quad+(t^2-t)z_{k_1+k_2+1}z_{k_3}\cdots z_{k_n}.
\end{align*}
Then we get \eqref{Eq:y-z1-zn-ast} by induction hypothesis.\qed

%%%%%%%%%%%%%%%%%%%%%%%%%%%%%%%%%%%%%%%%%%%%%%%%%%%%%%%%%%%%%%%%%%%%%%%

\section{A generating function for sums of $t$-MZVs}\label{Sec:Sum}

In this section, we study sums of $t$-MZVs of fixed weight, depth and height.

\subsection{A generating function for sums of $t$-MZVs}

As in Section \ref{Sec:Intro}, for positive integers $k,n,s$ with $k\geqslant n+s$ and $n\geqslant s$, let $I_0(k,n,s)$ be the set of all admissible indices of weight $k$, depth $n$ and height $s$, and let
$$X_0(k,n,s)=\sum\limits_{\mathbf{k}\in I_0(k,n,s)}\zeta^t(\mathbf{k})$$
be the sum of $t$-MZVs of weight $k$, depth $n$ and height $s$. We define a generating function of these sums
$$\Phi_0(u,v,w)=\sum\limits_{k\geqslant n+s,n\geqslant s\geqslant 1}X_0(k,n,s)u^{k-n-s}v^{n-s}w^{2s-2},$$
where $u,v,w$ are variables. To represent this function via hypergeometric functions, we follow the procedures given in \cite{Aoki-Kombu-Ohno,Aoki-Ohno} with more details.

For a sequence $\mathbf{k}=(k_1,k_2,\ldots,k_n)$ of positive integers, we define
$$\Li_{\mathbf{k}}(t,z)=\sum\limits_{\mathbf{p}} t^{n-\dep(\mathbf{p})}\Li_{\mathbf{p}}(z),$$
where $\mathbf{p}$ runs over all sequences of the form
$$\mathbf{p}=(k_1\Box k_2\Box\cdots\Box k_n)$$
in which each $\Box$ is filled by the comma, or the plus $+$, and
$$\Li_\mathbf{k}(z)=\sum\limits_{m_1>m_2>\cdots>m_n>0}\frac{z^{m_1}}{m_1^{k_1}m_2^{k_2}\cdots m_n^{k_n}}$$
is the multiple polylogarithm. Then the right-hand side of the definition equation of $\Li_{\mathbf{k}}(t,z)$ converges locally uniformly in the domain $|z|<1$, and converges in the point $z=1$ if $\mathbf{k}$ is admissible with the value of it coincides with $\zeta^{t}(\mathbf{k})$. By the iterated integral representation
$$\Li_{\mathbf{k}}(t,z)=\int\limits_{z>z_1>\cdots>z_k>0}\prod\limits_{i=1}^kf_i(z_i)dz_i$$
with
$$f_i(z)=\begin{cases}
\frac{t}{z}+\frac{1}{1-z} & \text{if\;} i=k_1,k_1+k_2,\ldots,k_1+\cdots+k_{n-1},\\
\frac{1}{1-z} & \text{if\;} i=k=\wt(\mathbf{k}),\\
\frac{1}{z} & \text{otherwise},
\end{cases}$$
we get
\begin{align}
\frac{d}{dz}\Li_{k_1,k_2,\ldots,k_n}(t,z)=\begin{cases}
\frac{1}{z}\Li_{k_1-1,k_2,\ldots,k_n}(t,z) & \text{if\;} k_1\geqslant 2,\\
\left(\frac{t}{z}+\frac{1}{1-z}\right)\Li_{k_2,\ldots,k_n}(t,z) & \text{if\;} n\geqslant 2, k_1=1, \\
\frac{1}{1-z} & \text{if\;} n=k_1=1.
\end{cases}
\label{Eq:Diff-t-MPL}
\end{align}

Now for nonnegative integers $k,n,s$, we denote by $I(k,n,s)$ the set of all indices of weight $k$, depth $n$ and height $s$, and define sums
\begin{align*}
X(k,n,s;z)=\sum\limits_{\mathbf{k}\in I(k,n,s)}\Li_{\mathbf{k}}(t,z),\\
X_0(k,n,s;z)=\sum\limits_{\mathbf{k}\in I_0(k,n,s)}\Li_{\mathbf{k}}(t,z).
\end{align*}
If the index set $I(k,n,s)$ or $I_0(k,n,s)$ is an empty set, the corresponding sum is treated as zero, and we use the convention $X(0,0,0;z)=1$. Using \eqref{Eq:Diff-t-MPL}, we find for any integers $k,n,s$ with $k\geqslant n+s$, the following equations hold
\begin{align}
&\frac{d}{dz}X_0(k,n,s;z)=\frac{1}{z}\left(X(k-1,n,s-1;z)-X_0(k-1,n,s-1;z)\right.\nonumber\\
&\qquad\qquad\left.+X_0(k-1,n,s,;z)\right)\qquad\qquad (n\geqslant s\geqslant 1),
\label{Eq:Diff-X0}\\
&\frac{d}{dz}\left(X(k,n,s;z)-X_0(k,n,s;z)\right)=\left(\frac{t}{z}+\frac{1}{1-z}\right)X(k-1,n-1,s;z)\nonumber\\
&\qquad\qquad\qquad\qquad\qquad\qquad\qquad\qquad\quad (n\geqslant s\geqslant 0,n\geqslant 2).
\label{Eq:Diff-X-X0}
\end{align}

Then we define the following generating functions
\begin{align*}
&\Phi(z)=\Phi(u,v,w;z)=\sum\limits_{k,n,s\geqslant 0}X(k,n,s;z)u^{k-n-s}v^{n-s}w^{2s},\\
&\Phi_0(z)=\Phi_0(u,v,w;z)=\sum\limits_{k,n,s\geqslant 0}X_0(k,n,s;z)u^{k-n-s}v^{n-s}w^{2s-2}.
\end{align*}
Using \eqref{Eq:Diff-X0} and \eqref{Eq:Diff-X-X0}, we get
\begin{align*}
&\frac{d}{dz}\Phi_0(z)=\frac{1}{vz}\left(\Phi(z)-1-w^2\Phi_0(z)\right)+\frac{u}{z}\Phi_0(z),\\
&\frac{d}{dz}\left(\Phi(z)-w^2\Phi_0(z)\right)=\left(\frac{t}{z}+\frac{1}{1-z}\right)v\left(\Phi(z)-1\right)+\frac{v}{1-z}.
\end{align*}
Eliminating $\Phi(z)$, we get the differential equation that satisfied by $\Phi_0(z)$.

\begin{prop}\label{Prop:Diff-Equation}
The formal power series $\Phi_0(z)$ is a unique power series solution vanishing at $z=0$ of the differential equation
\begin{align}
z^2(1-z)\frac{d^{2}f}{dz^2}+z\{(1-u-vt)(1-z)-vz\}\frac{df}{dz}+[t+(1-t)z](uv-w^2)f=z.
\label{Eq:Diff-Equation}
\end{align}
Hence it converges locally uniformly and defines a holomorphic function in $|z|<1$.
\end{prop}

\proof The homogeneous equation of \eqref{Eq:Diff-Equation} is a second-order linear differential equation of Fuchsian type with singularities at $z=0,1,\infty$. The characteristic equation at $z=0$ is
\begin{align}
\lambda^2-(u+vt)\lambda+t(uv-w^2)=0,
\label{Eq:ChaEqu-0}
\end{align}
which induces that the characteristic exponents at $z=0$ are not positive integers for generic $u,v,w$. Hence \eqref{Eq:Diff-Equation} has a unique formal power series solution vanishing at $z=0$ and the solution should converge locally in $|z|<1$.
\qed

Let $\alpha,\beta$ be roots of the characteristic equation \eqref{Eq:ChaEqu-0}. Then we have
$$\alpha+\beta=u+vt,\qquad \alpha\beta=t(uv-w^2).$$
The characteristic equation of the homogeneous equation of \eqref{Eq:Diff-Equation} at $z=\infty$ is
\begin{align}
\lambda^2+(u+(t-1)v)\lambda+(t-1)(uv-w^2)=0.
\label{Eq:ChaEqu-infty}
\end{align}
Let $\gamma_1,\gamma_2$ be roots of \eqref{Eq:ChaEqu-infty}. Then we get
$$\gamma_1+\gamma_2=-u-(t-1)v,\qquad \gamma_1\gamma_2=(t-1)(uv-w^2).$$
Note that we have
$$\alpha+\beta+\gamma_1+\gamma_2=v.$$

Assume that $f=\sum\limits_{n=1}^\infty a_nz^n$ is a solution of \eqref{Eq:Diff-Equation}, then we have
\begin{align*}
&a_1=\frac{1}{(1-u)(1-vt)-tw^2},\\
&a_n=\frac{(n-1)(n-1-u+(1-t)v)-(1-t)(uv-w^2)}{(n-u)(n-vt)-tw^2}a_{n-1}\quad (n\geqslant 2).
\end{align*}
Therefore we obtain the following result, which is just \cite[Proposition 2.2]{Aoki-Kombu-Ohno} when $t=1$.

\begin{prop}\label{Prop:Solu-PowerSeries}
We have
$$\Phi_0(u,v,w)=\sum\limits_{n=1}^\infty a_n,$$
where
$$a_n=\frac{\prod\limits_{m=1}^{n-1}(m+\gamma_1)(m+\gamma_2)}{\prod\limits_{m=1}^n(m-\alpha)(m-\beta)}.$$
\end{prop}

On the other hand, the characteristic equation of the homogeneous equation of \eqref{Eq:Diff-Equation} at $z=1$ is
$$\lambda^2+(v-1)\lambda=0,$$
which has roots $0$ and $1-v$. Then using Riemann's $P$-function, the system of fundamental solutions of the homogeneous equation of \eqref{Eq:Diff-Equation} is
\begin{align*}
&G_1(z)=z^{\alpha}F\left({\gamma_1+\alpha,\gamma_2+\alpha\atop\alpha-\beta+1};z\right),\\
&G_2(z)=z^{\beta}F\left({\gamma_1+\beta,\gamma_2+\beta\atop \beta-\alpha+1};z\right),
\end{align*}
where $F\left({a,b\atop c};z\right)$ is the Gauss hypergeometric function defined as
$$F\left({a,b\atop c};z\right)=\sum\limits_{n=0}^\infty \frac{(a)_n(b)_n}{n!(c)_n}z^n,$$
with the Pochhammer symbol $(a)_n$ given by
$$(a)_n=\begin{cases}
1& \text{if\;} n=0,\\
a(a+1)\cdots (a+n-1) & \text{if\;} n>0.
\end{cases}$$

We want to employ the method of variation of constants to get a solution of \eqref{Eq:Diff-Equation}. Then it needs to compute the Wronskian of $G_1(z)$ and $G_2(z)$.

\begin{lem}\label{Lem:Wronskian}
The Wronskian of $G_1(z)$ and $G_2(z)$ is
$$\begin{vmatrix}
G_1(z) & G_2(z) \\
G_1'(z) & G_2'(z)
\end{vmatrix}=(\beta-\alpha)z^{\alpha+\beta-1}(1-z)^{-v}.$$
\end{lem}

\proof Denote the Wronskian by $W(z)$, then it is direct to get
$$W(z)=(\beta-\alpha)z^{\alpha+\beta-1}F_1(z)F_2(z)+z^{\alpha+\beta}(F_1(z)F_2'(z)-F_1'(z)F_2(z)),$$
where
\begin{align*}
F_1(z)=F\left({\gamma_1+\alpha,\gamma_2+\alpha\atop \alpha-\beta+1};z\right),\quad F_2(z)=F\left({\gamma_1+\beta,\gamma_2+\beta\atop\beta-\alpha+1};z\right).
\end{align*}
Using the formula
$$F\left({a,b\atop c};z\right)+\frac{z}{a}\frac{d}{dz}F\left({a,b\atop c};z\right)=F\left({a+1,b\atop c};z\right),$$
we get
\begin{align*}
&(\beta-\alpha)F_1(z)F_2(z)+z(F_1(z)F_2'(z)-F_1'(z)F_2(z))\\
=&(\gamma_1+\beta)F\left({\gamma_1+\alpha,\gamma_2+\alpha\atop\alpha-\beta+1};z\right)
F\left({\gamma_1+\beta+1,\gamma_2+\beta\atop\beta-\alpha+1};z\right)\\
&-(\gamma_1+\alpha)F\left({\gamma_1+\alpha+1,\gamma_2+\alpha\atop\alpha-\beta+1};z\right)
F\left({\gamma_1+\beta,\gamma_2+\beta\atop\beta-\alpha+1};z\right).
\end{align*}
Applying the formula (see \cite{Bailey} for example)
\begin{align}
F\left({a,b\atop c};z\right)=(1-z)^{c-a-b}F\left({c-a,c-b\atop c};z\right),
\label{Eq:GaussHyp-c-ab}
\end{align}
we get
\begin{align*}
W(z)=&z^{\alpha+\beta-1}(1-z)^{-v}\\
&\times \left[(\gamma_1+\beta)F\left({\gamma_1+\alpha,\gamma_2+\alpha\atop \alpha-\beta+1};z\right)F\left({-\gamma_1-\alpha,1-\gamma_2-\alpha\atop\beta-\alpha+1};z\right)\right.\\
&\left.-(\gamma_1+\alpha)F\left({\gamma_1+\beta,\gamma_2+\beta\atop\beta-\alpha+1};z\right)
F\left({-\gamma_1-\beta,1-\gamma_2-\beta\atop\alpha-\beta+1};z\right)\right].
\end{align*}
Now the result follows from
\begin{align}
&(\gamma_1+\beta)F\left({\gamma_1+\alpha,\gamma_2+\alpha\atop\alpha-\beta+1};z\right)F\left({-\gamma_1-\alpha,1-\gamma_2-\alpha\atop \beta-\alpha+1};z\right)\nonumber\\
&\quad-(\gamma_1+\alpha)F\left({\gamma_1+\beta,\gamma_2+\beta\atop\beta-\alpha+1};z\right)
F\left({-\gamma_1-\beta,1-\gamma_2-\beta\atop\alpha-\beta+1};z\right)=\beta-\alpha.
\label{Eq:Gauss-Hyper-Duality}
\end{align}
A proof of \eqref{Eq:Gauss-Hyper-Duality} will be given in Appendix \ref{AppSec:Proof-Gauss}.
\qed

Hence by the method of variation of constants, we find that
$$\Phi_0(u,v,w;z)=C_1(z)G_1(z)+C_2(z)G_2(z),$$
where
\begin{align*}
C_1(z)=&\frac{1}{\alpha-\beta}\int_0^zx^{-\alpha}(1-x)^{v-1}F\left({\gamma_1+\beta,\gamma_2+\beta\atop\beta-\alpha+1};x\right)dx,\\
C_2(z)=&\frac{1}{\beta-\alpha}\int_0^zx^{-\beta}(1-x)^{v-1}F\left({\gamma_1+\alpha,\gamma_2+\alpha\atop\alpha-\beta+1};x\right)dx.
\end{align*}
Here and below we assume that $\Re(\alpha)<1$ and $\Re(\beta)<1$. Let $z=1$, and use Gaussian summation formula (see \cite{Bailey} for example)
$$F\left({a,b\atop c};1\right)=\frac{\Gamma(c)\Gamma(c-a-b)}{\Gamma(c-a)\Gamma(c-b)}\qquad (\Re(c-a-b)>0),$$
we get the following result.

\begin{prop}\label{Prop:Phi-GaussHyp}
We have
\begin{align*}
&\Phi_0(u,v,w)\\
=&\frac{\Gamma(\alpha-\beta)\Gamma(1-v)}{\Gamma(1-\gamma_1-\beta)\Gamma(1-\gamma_2-\beta)}\int_0^1z^{-\alpha}(1-z)^{v-1}
F\left({\gamma_1+\beta,\gamma_2+\beta\atop\beta-\alpha+1};z\right)dz\\
&+\frac{\Gamma(\beta-\alpha)\Gamma(1-v)}{\Gamma(1-\gamma_1-\alpha)\Gamma(1-\gamma_2-\alpha)}\int_0^1z^{-\beta}(1-z)^{v-1}
F\left({\gamma_1+\alpha,\gamma_2+\alpha\atop\alpha-\beta+1};z\right)dz.
\end{align*}
\end{prop}

Recall from \cite{Bailey} that the generalized hypergeometric function $_3F_2$ is defined by
$$\;_3F_2\left({a_1,a_2,a_3\atop b_1,b_2};z\right)=\sum\limits_{n=0}^\infty \frac{(a_1)_n(a_2)_n(a_3)_n}{n!(b_1)_n(b_2)_n}z^n.$$
Then using some connection formula of Gauss hypergeometric functions, we get the following result.

\begin{thm}\label{Thm:Sum}
We have
\begin{align*}
\Phi_0(u,v,w)=&\frac{1}{1-v}\int_0^1(1-z)^{-\beta}F\left({1-\gamma_1-\beta, 1-\gamma_2-\beta\atop 2-v};z\right)dz\\
=&\frac{1}{(1-v)(1-\beta)}\;_3F_2\left({1-\gamma_1-\beta,1-\gamma_2-\beta,1\atop 2-v,2-\beta};1\right),
\end{align*}
where $\alpha,\beta$ are determined by $\alpha+\beta=u+vt$, $\alpha\beta=t(uv-w^2)$, and $\gamma_1,\gamma_2$ are determined by $\gamma_1+\gamma_2=-u-(t-1)v$, $\gamma_1\gamma_2=(t-1)(uv-w^2)$.
\end{thm}

\proof Using the formula \eqref{Eq:GaussHyp-c-ab}, we can rewrite the formula for $\Phi_0(u,v,w)$ in Proposition \ref{Prop:Phi-GaussHyp} as
\begin{align*}
&\Phi_0(u,v,w)=\int_0^1dz\cdot z^{-\beta}\\
&\quad \times\left[\frac{\Gamma(\beta-\alpha)\Gamma(1-v)}{\Gamma(1-\gamma_1-\alpha)\Gamma(1-\gamma_2-\alpha)}
F\left({1-\gamma_1-\beta,1-\gamma_2-\beta\atop\alpha-\beta+1};z\right)\right.\\
&\quad\left.+\frac{\Gamma(\alpha-\beta)\Gamma(1-v)}{\Gamma(1-\gamma_1-\beta)\Gamma(1-\gamma_2-\beta)}z^{\beta-\alpha}
F\left({1-\gamma_1-\alpha,1-\gamma_2-\alpha\atop\beta-\alpha+1};z\right)\right].
\end{align*}
Then applying the connection formula (see \cite{Bailey} for example)
\begin{align*}
F\left({a,b\atop c};z\right)=&\frac{\Gamma(c)\Gamma(c-a-b)}{\Gamma(c-a)\Gamma(c-b)}F\left({a,b\atop 1+a+b-c};1-z\right)\\
&+\frac{\Gamma(c)\Gamma(a+b-c)}{\Gamma(a)\Gamma(b)}(1-z)^{c-a-b}F\left({c-a,c-b\atop 1+c-a-b};1-z\right)
\end{align*}
with $a=1-\gamma_1-\beta$, $b=1-\gamma_2-\beta$ and $c=2-v$, we get
$$\Phi_0(u,v,w)=\frac{1}{1-v}\int_0^1 z^{-\beta}F\left({1-\gamma_1-\beta,1-\gamma_2-\beta\atop 2-v};1-z\right)dz.$$
Now it is easy to get the result.
\qed

Let $t=0$ in Theorem \ref{Thm:Sum}, we get Ohno-Zagier relation by applying Gaussian summation formula. And let $t=1$ in Theorem \ref{Thm:Sum}, we get \cite[Proposition 3.1]{Aoki-Kombu-Ohno}.

\subsection{The case $uv=w^2$}

Setting $uv=w^2$, then
$$\alpha+\beta=u+vt,\quad \alpha\beta=0$$
and
$$\gamma_1+\gamma_2=-u-(t-1)v,\quad \gamma_1\gamma_2=0.$$
Hence we take $\alpha=0$, $\beta=u+vt$, $\gamma_1=0$ and $\gamma_2=-u-(t-1)v$. Then by Theorem \ref{Thm:Sum}, we  get
\begin{align*}
\Phi_0(u,v,w)|_{uv=w^2}=&\frac{1}{(1-v)(1-\beta)}\;_3F_2\left({1-\beta,1-\gamma_2-\beta,1\atop 2-v,2-\beta};1\right)\\
=&\frac{1}{(1-v)(1-\beta)}\sum\limits_{n=0}^\infty \frac{(1-v)_n(1-\beta)_n}{(2-v)_n(2-\beta)_n}\\
=&\sum\limits_{n=1}^\infty\frac{1}{(n-\beta)(n-v)}
=\sum\limits_{n=1}^\infty\sum\limits_{m=1}^\infty\sum\limits_{l=1}^\infty\frac{1}{n^{m+l}}\beta^{m-1}v^{l-1}\\
=&\sum\limits_{m,l=1}^\infty\zeta(m+l)\sum\limits_{i=0}^{m-1}\binom{m-1}{i}u^{m-1-i}v^{i+l-1}t^{i}\\
=&\sum\limits_{k\geqslant n+1,n\geqslant 1}\left(\sum\limits_{i=0}^{n-1}\binom{k-n+i-1}{i}t^i\right)\zeta(k)u^{k-n-1}v^{n-1}.
\end{align*}
Therefore we have
\begin{align}
\sum\limits_{\wt(\mathbf{k})=k,\dep(\mathbf{k})=n\atop \mathbf{k}:\text{admissible}}\zeta^t(\mathbf{k})=\left(\sum\limits_{i=0}^{n-1}\binom{k-n+i-1}{i}t^i\right)\zeta(k).
\label{Eq:Sum-Formula-New}
\end{align}
Note that \eqref{Eq:Sum-Formula-New} is just \eqref{Eq:Sum-Formula}, since one can check that
$$\sum\limits_{i=0}^{n-1}\binom{k-1}{i}t^i(1-t)^{n-1-i}=\sum\limits_{i=0}^{n-1}\binom{k-n+i-1}{i}t^i\qquad (n\in\mathbb{N}).$$
Therefore we give a new proof of the sum formula of $t$-MZVs.

\subsection{The case $v=0$}

Setting $v=0$, then by Theorem \ref{Thm:Sum}, we have
\begin{align*}
\Phi_0(u,v,w)|_{v=0}=&\frac{1}{1-\beta}\;_3F_2\left({1-\gamma_1-\beta,1-\gamma_2-\beta,1\atop 2,2-\beta};1\right)\\
=&\frac{1}{1-\beta}\sum\limits_{n=0}^\infty \frac{(1-\gamma_1-\beta)_n(1-\gamma_2-\beta)_n}{(n+1)!(2-\beta)_n}\\
=&\frac{1}{(\gamma_1+\beta)(\gamma_2+\beta)}\sum\limits_{n=1}^\infty\frac{(-\gamma_1-\beta)_n(-\gamma_2-\beta)_n}{n!(1-\beta)_n}\\
=&\frac{1}{(\gamma_1+\beta)(\gamma_2+\beta)}\left\{F\left({-\gamma_1-\beta,-\gamma_2-\beta\atop 1-\beta};1\right)-1\right\}.
\end{align*}
Using Gaussian summation formula, we get
$$\Phi_0(u,v,w)|_{v=0}=\frac{1}{(\gamma_1+\beta)(\gamma_2+\beta)}
\left\{\frac{\Gamma(1-\alpha)\Gamma(1-\beta)}{\Gamma(1+\gamma_1)(1+\gamma_2)}-1\right\}.$$
Applying the expansion
$$\Gamma(1-x)=\exp\left(\gamma x+\sum\limits_{n=2}^\infty \frac{\zeta(n)}{n}x^n\right)$$
with $\gamma$ the Euler-Mascheroni constant, we find the sum $X_0(k,n,n)$ can be expressed as a polynomials of Riemann zeta values with $\mathbb{Q}[t]$ coefficients as follows.

\begin{cor}\label{Cor:v-0}
Let $\alpha,\beta,\delta_1$ and $\delta_2$ be determined by
$$\alpha+\beta=u,\quad \alpha\beta=-tw^2$$
and
$$\delta_1+\delta_2=u,\quad \delta_1\delta_2=(1-t)w^2.$$
We have
\begin{align*}
\sum\limits_{k\geqslant 2s, s\geqslant 1}X_0(k,n,n)u^{k-2s}w^{2s}
=\exp\left[\sum\limits_{n=2}^\infty\frac{\zeta(n)}{n}
(\alpha^n+\beta^n-\delta_1^n-\delta_2^n)\right]-1.
\end{align*}
\end{cor}

Set $t=1$ in Corollary \ref{Cor:v-0}, we get \cite[Theorem 4.2]{Aoki-Kombu-Ohno}.

\subsection{The case $w=0$}

Finally, we set $w=0$. Then we have
$$\alpha+\beta=u+vt,\quad \alpha\beta=tuv$$
and
$$\gamma_1+\gamma_2=-u-(t-1)v,\quad \gamma_1\gamma_2=(t-1)uv.$$
Hence we take $\alpha=u$, $\beta=vt$, $\gamma_1=-u$ and $\gamma_2=(1-t)v$. By Proposition \ref{Prop:Solu-PowerSeries}, we have
\begin{align*}
\Phi_0(u,v,w)|_{w=0}=&\sum\limits_{n=1}^\infty\frac{\prod\limits_{m=1}^{n-1}(m-u)(m+(1-t)v)}{\prod\limits_{m=1}^n(m-u)(m-tv)}\\
=&\sum\limits_{n=1}^\infty\frac{\prod\limits_{m=1}^{n-1}(m+(1-t)v)}{(n-u)\prod\limits_{m=1}^n(m-tv)}.
\end{align*}

When $t=0$, we have
\begin{align*}
\Phi_0(u,v,w)|_{w=0}=&-\frac{1}{uv}\sum\limits_{n=1}^\infty\frac{(-u)_n(v)_n}{n!(1-u)_n}=\frac{1}{uv}\left\{1-F\left({-u,v\atop 1-u};1\right)\right\}\\
=&\frac{1}{uv}\left(1-\frac{\Gamma(1-u)\Gamma(1-v)}{\Gamma(1)\Gamma(1-u-v)}\right)\\
=&\frac{1}{uv}\left\{1-\exp\left[\sum\limits_{n=2}^\infty\frac{\zeta(n)}{n}\left(u^n+v^n-(u+v)^n\right)\right]\right\}.
\end{align*}
Then we have
\begin{align*}
&\sum\limits_{k\geqslant n+1,n\geqslant 1}\zeta(k-n+1,\underbrace{1,\ldots,1}_{n-1})u^{k-n}v^n\\
=&1-\exp\left\{\sum\limits_{n=2}^\infty\frac{\zeta(n)}{n}\left[u^n+v^n-(u+v)^n\right]\right\},
\end{align*}
which is exactly the Aomoto-Drinfel'd-Zagier relation.

When $t\neq 0$, we obtain a formula for height one $t$-MZVs as follows, which recovers \cite[Proposition 4.3]{Aoki-Kombu-Ohno} in the case $t=1$.

\begin{cor}\label{Cor:w-0}
For any positive integers $i,j$ and a nonzero variable $t$, we have
$$\zeta^t(i+1,\underbrace{1,\ldots,1}_{j-1})=\sum\limits_{n\geqslant m\geqslant 1}(-1)^{m-1}\binom{n-1}{m-1}\prod\limits_{l=1}^{n-1}\left(1+\frac{1-t}{t}\frac{m}{l}\right)\frac{t^{j-1}}{n^im^j}.$$
\end{cor}

\proof If  $t\neq 0$, we get
\begin{align*}
&\Phi_0(u,v,w)|_{w=0}=\sum\limits_{n=1}^\infty\frac{\prod\limits_{m=1}^{n-1}\left(m+\frac{1-t}{t}tv\right)}{(n-u)\prod\limits_{m=1}^n(m-tv)}\\
=&\sum\limits_{n=1}^\infty\sum\limits_{m=1}^n\frac{(-1)^{m-1}}{(n-u)(m-tv)}\frac{\prod\limits_{l=1}^{n-1}\left(l+\frac{1-t}{t}m\right)}{(m-1)!(n-m)!}\\
=&\sum\limits_{i,j=1}^\infty\left\{\sum\limits_{n\geqslant m\geqslant 1}(-1)^{m-1}\binom{n-1}{m-1}\prod\limits_{l=1}^{n-1}\left(1+\frac{1-t}{t}\frac{m}{l}\right)\frac{t^{j-1}}{n^im^j}\right\}u^{i-1}v^{j-1}.
\end{align*}
Then we get the result.
\qed

\subsection{Further remarks}

In \cite{Ohno-Zagier}, setting $u=-v$ in Ohno-Zagier relation, the authors obtained the Le-Murakami relation, which was first proved in \cite{Le-Murakami}. In \cite{Aoki-Kombu-Ohno,Aoki-Ohno}, setting $u=v$ in the formula of the generating function, the authors obtained a formula for sums of MZSVs of fixed weight and height. Then it is natural to ask whether analogous results for sums of $t$-MZVs of fixed weight and height can be deduced from Proposition \ref{Prop:Solu-PowerSeries} and Theorem \ref{Thm:Sum} or not.

Also using some formulas of hypergeometric functions, one can obtain some duality for MZSVs (See \cite{Kaneko-Ohno,Li2012,Yamazaki,Yamazaki-thesis}). More precisely, when $t=1$, in \cite{Li2012} or in \cite{Yamazaki-thesis}, the authors proved that the difference
$$u\Phi_0(-u,v,w)-v\Phi_0(-v,u,w)$$
can be expressed by gamma functions.  Then one may ask whether a similar duality for $t$-MZVs can be deduced from Theorem \ref{Thm:Sum} or not.

%%%%%%%%%%%%%%%%%%%%%%%%%%%%%%%%%%%%%%%%%%%%%%%%%%%%%%%%%%%%%%%%%%%%%%%%%%%%%
%%%%%%%%%%%%%%%%%%%%%%%%%%%%%%%%%%%%%%%%%%%%%%%%%%%%%%%%%%%%%%%%%%%%%%%%%%%%%

\appendix

%%%%%%%%%%%%%%%%%%%%%%%%%%%%%%%%%%%%%%%%%%%%%%%%%%%%%%%%%%%%%%%%%%%%%%%%%%%%%
%%%%%%%%%%%%%%%%%%%%%%%%%%%%%%%%%%%%%%%%%%%%%%%%%%%%%%%%%%%%%%%%%%%%%%%%%%%%%

\section{Another proof of the sum formula}\label{AppSec:Proof-Sum}

In this section, we derive the sum formula \eqref{Eq:Sum-Formula-New} from that of MZVs. We prepare a lemma.

\begin{lem}\label{Lem:Sigma}
For any $w_1,w_2\in\mathfrak{h}_t$, we have
$$\sigma_t(w_1\shuffle w_2)=\sigma_t(w_1)\shuffle \sigma_t(w_2).$$
\end{lem}

\proof We may assume that $w_1,w_2\in A^{\ast}$. Then it is easy to get the result by induction on $|w_1|+|w_2|$.\qed

Now let $k,n$ be integers with $k>n\geqslant 1$. By the sum formula of MZVs, we have
$$x(x^{k-n-1}\shuffle y^{n-1})y-x^{k-1}y\in \ker Z.$$
While on $\mathfrak{h}_t^0$ we have $Z=Z_t\circ S_{-t}$. Hence we get
$$S_{-t}(x(x^{k-n-1}\shuffle y^{n-1})y)-S_{-t}(x^{k-1}y)\in \ker Z_t.$$
Since
\begin{align*}
&S_{-t}(x(x^{k-n-1}\shuffle y^{n-1})y)=x\sigma_{-t}(x^{k-n-1}\shuffle y^{n-1})y\\
=&x(x^{k-n-1}\shuffle (-tx+y)^{n-1})y\qquad (\text{by Lemma \ref{Lem:Sigma}})\\
=&\sum\limits_{i=0}^{n-1}(-t)^ix(x^{k-n-1}\shuffle x^i\shuffle y^{n-1-i})y\\
=&\sum\limits_{i=0}^{n-1}(-t)^i\binom{k-n+i-1}{i}x(x^{k-n+i-1}\shuffle y^{n-1-i})y\\
=&\sum\limits_{i=1}^n(-t)^{n-i}\binom{k-i-1}{n-i}x(x^{k-i-1}\shuffle y^{i-1})y,
\end{align*}
we get
$$\sum\limits_{i=1}^n(-t)^{n-i}\binom{k-i-1}{n-i}G_0(k,i)=\zeta(k),$$
where
$$G_0(k,i)=\sum\limits_{\wt(\mathbf{k})=k,\dep(\mathbf{k})=i\atop \mathbf{k}:\text{admissible}}\zeta^t(\mathbf{k}).$$
We rewrite the above equation as
$$\sum\limits_{i=1}^n\binom{n}{i}(-t)^{-i}(k-1-i)!i!G_0(k,i)=n!(k-n-1)!(-t)^{-n}\zeta(k).$$
Then using the binomial inversion formula
$$a_n=\sum\limits_{i=1}^n\binom{n}{i}b_i \Longleftrightarrow b_n=\sum\limits_{i=1}^n(-1)^{n-i}\binom{n}{i}a_i,$$
we get
$$(-t)^{-n}(k-1-n)!n!G_0(k,n)=\sum\limits_{i=1}^n(-1)^{n-i}\binom{n}{i}i!(k-i-1)!(-t)^{-i}\zeta(k).$$
Rewriting the above equation, we get
$$G_0(k,n)=\sum\limits_{i=1}^n\binom{k-1-i}{n-i}t^{n-i}\zeta(k).$$
which is just \eqref{Eq:Sum-Formula-New}.

%%%%%%%%%%%%%%%%%%%%%%%%%%%%%%%%%%%%%%%%%%%%%%%%%%%%%%%%%%%%%%%%%%%%%%%%%%%%%%%%%%%%%%%%%%%%%%%%%%%%%%%%%%%%%%%%%%

\section{A proof of \eqref{Eq:Gauss-Hyper-Duality}}\label{AppSec:Proof-Gauss}

In this section, we give a proof of \eqref{Eq:Gauss-Hyper-Duality}. It is enough to show for any positive integer $n$, it holds that
\begin{align*}
&\sum\limits_{i+j=n\atop i,j\geqslant 0}(\gamma_1+\beta)\frac{(\gamma_1+\alpha)_i(\gamma_2+\alpha)_i(-\gamma_1-\alpha)_j(1-\gamma_2-\alpha)_j}{i!j!(\alpha-\beta+1)_i(\beta-\alpha+1)_j}\\
=&\sum\limits_{i+j=n\atop i,j\geqslant 0}(\gamma_1+\alpha)\frac{(\gamma_1+\beta)_i(\gamma_2+\beta)_i(-\gamma_1-\beta)_j(1-\gamma_2-\beta)_j}{i!j!(\beta-\alpha+1)_i(\alpha-\beta+1)_j}.
\end{align*}
Then it is sufficient to show that
\begin{align*}
&\sum\limits_{i+j=n\atop i,j\geqslant 0}(-1)^j\frac{(\gamma_1+\alpha-j+1)_{n-1}(\gamma_2+\alpha-j)_n}{i!j!(\alpha-\beta-j)_{n+1}}\\
=&\sum\limits_{i+j=n\atop i,j\geqslant 0}(-1)^j\frac{(\gamma_1+\beta-i+1)_{n-1}(\gamma_2+\beta-i)_n}{i!j!(\alpha-\beta-j)_{n+1}}.
\end{align*}
Assume that $\alpha-\beta=x$, hence it is equivalent to show
\begin{align*}
&\sum\limits_{i+j=n\atop i,j\geqslant 0}(-1)^j\frac{(\gamma_1+\alpha-j+1)_{n-1}(\gamma_2+\alpha-j)_n}{i!j!(x-j)_{n+1}}\\
=&\sum\limits_{i+j=n\atop i,j\geqslant 0}(-1)^j\frac{(\gamma_1+\alpha-x-i+1)_{n-1}(\gamma_2+\alpha-x-i)_n}{i!j!(x-j)_{n+1}}.
\end{align*}
Finally, we get \eqref{Eq:Gauss-Hyper-Duality} from the following lemma by setting $\gamma_1+\alpha=y_1$ and $\gamma_2+\alpha=y_2$.

\begin{lem}
For any positive integer $n$ and any variables $x,y_1,y_2$, we have
\begin{align}
&\sum\limits_{j=0}^n(-1)^j\frac{(y_1-j+1)_{n-1}(y_2-j)_n}{j!(n-j)!(x-j)_{n+1}}\nonumber\\
=&\sum\limits_{j=0}^n(-1)^j\frac{(y_1-x-n+j+1)_{n-1}(y_2-x-n+j)_n}{j!(n-j)!(x-j)_{n+1}}.
\label{Eq:Gauss-Hyper-Explicit}
\end{align}
\end{lem}

\proof
It is easy to check that \eqref{Eq:Gauss-Hyper-Explicit} holds for $n=1$. Now we assume that $n\geqslant 2$.

Set
\begin{align*}
f(x)=&\sum\limits_{j=0}^n(-1)^j\frac{(y_1-x-n+j+1)_{n-1}(y_2-x-n+j)_n}{j!(n-j)!(x-j)_{n+1}}\\
&-\sum\limits_{j=0}^n(-1)^j\frac{(y_1-j+1)_{n-1}(y_2-j)_n}{j!(n-j)!(x-j)_{n+1}}\in\mathbb{Q}[y_1,y_2](x).
\end{align*}
We show that $f(x)$ is a polynomial in $x$ with degree at most $n-2$. In fact, let $q_j(x)$ and $r_j(x)$ be polynomials with $\deg q_j(x)=n-2$, $\deg r_j(x)\leqslant n$ and
$$(y_1-x-n+j+1)_{n-1}(y_2-x-n+j)_n=(x-j)_{n+1}q_j(x)+r_j(x).$$
Then let $x=l$ with $l=j,j-1,\ldots,j-n$ in above equation, we get
$$r_j(l)=(y_1-l-n+j+1)_{n-1}(y_2-l-n+j)_n.$$
Hence we find
\begin{align*}
r_j(x)=&\sum\limits_{l=j-n}^j\frac{(-1)^{j-l}(y_1-l-n+j+1)_{n-1}(y_2-l-n+j)_n(x-j)_{n+1}}{(j-l)!(l-j+n)!(x-l)}\\
=&(x-j)_{n+1}\sum\limits_{l=0}^n(-1)^l\frac{(y_1+l-n+1)_{n-1}(y_2+l-n)_n}{l!(n-l)!(x+l-j)}.
\end{align*}
Since
$$\frac{1}{(x-j)_{n+1}}=\sum\limits_{l=0}^n(-1)^l\frac{1}{l!(n-l)!(x+l-j)},$$
we have
\begin{align*}
f(x)=&\sum\limits_{j=0}^n(-1)^j\frac{q_j(x)}{j!(n-j)!}+\sum\limits_{j,l=0}^n(-1)^{j+l}\frac{(y_1+l-n+1)_{n-1}(y_2+l-n)_n}{j!(n-j)!l!(n-l)!(x+l-j)}\\
&-\sum\limits_{j,l=0}^n(-1)^{j+l}\frac{(y_1-j+1)_{n-1}(y_2-j)_n}{j!(n-j)!l!(n-l)!(x+l-j)}\\
=&\sum\limits_{j=0}^n(-1)^j\frac{q_j(x)}{j!(n-j)!}.
\end{align*}
Therefore $f(x)$ is a polynomial in $x$ with degree at most $n-2$.

Thus to prove \eqref{Eq:Gauss-Hyper-Explicit}, it is sufficient to show that
$$f(y_1+i)=0$$
for $i=1,2,\ldots,n-1$. In fact, for any $i$, we have
\begin{align*}
f(y_1+i)=&\sum\limits_{j=0}^n(-1)^j\frac{(-i-n+j+1)_{n-1}(y_2-y_1-i-n+j)_n}{j!(n-j)!(y_1+i-j)_{n+1}}\\
&-\sum\limits_{j=0}^n(-1)^j\frac{(y_1-j+1)_{n-1}(y_2-j)_n}{j!(n-j)!(y_1+i-j)_{n+1}}.
\end{align*}
Using the partial fraction expansions
$$\frac{(y_2-y_1-i-n+j)_n}{(y_1+i-j)_{n+1}}=\sum\limits_{l=0}^n(-1)^l\frac{(y_2+l-n)_n}{l!(n-l)!(y_1+i-j+l)}$$
and
$$\frac{(y_1-j+1)_{n-1}}{(y_1+i-j)_{n+1}}=\sum\limits_{l=0}^n(-1)^l\frac{(-i-l+1)_{n-1}}{l!(n-l)!(y_1+i-j+l)},$$
we find
\begin{align*}
f(y_1+i)=&\sum\limits_{j,l=0}^n(-1)^{j+l}\frac{(-i-n+j+1)_{n-1}(y_2+l-n)_n}{j!(n-j)!l!(n-l)!(y_1+i-j+l)}\\
&-\sum\limits_{j,l=0}^n(-1)^{j+l}\frac{(-i-l+1)_{n-1}(y_2-j)_n}{j!(n-j)!l!(n-l)!(y_1+i-j+l)}\\
=&0,
\end{align*}
which finishes the proof.
\qed

%%%%%%%%%%%%%%%%%%%%%%%%%%%%%%%%%%%%%%%%%%%%%%%%%%%%%%%%%%%%%%%%%%%%%%%%%%%%%%%%%%%%%%%%%%%%%%%%%%%%%%%%%%%%%%%%%%%%

\end{document}